\documentclass[12pt]{amsart}

\newtheorem{theorem}{Theorem}[section]
\newtheorem*{maintheorem}{Theorem}

\newtheorem{lemma}[theorem]{Lemma}
\newtheorem{proposition}[theorem]{Proposition}
\newtheorem{corollary}[theorem]{Corollary}
\newtheorem*{maincorollary}{Corollary}

\theoremstyle{remark}

\newtheorem*{acknowlegments}{Acknowlegments}
\newtheorem{example}[theorem]{Example}

\theoremstyle{definition}
\newtheorem{definition}[theorem]{Definition}

\newcommand{\adj}{\mathrm{Adj}}
\newcommand{\et}{\quad\mbox{and}\quad}

\newcommand{\bC}{\mathbb{C}}
\newcommand{\bP}{\mathbb{P}}
\newcommand{\bQ}{\mathbb{Q}}
\newcommand{\bR}{\mathbb{R}}
\newcommand{\bZ}{\mathbb{Z}}
\newcommand{\cC}{{\mathcal{C}}}

\newcommand{\cM}{{\mathcal{M}}}

\newcommand{\cU}{{\mathcal{U}}}
\newcommand{\cV}{{\mathcal{V}}}
\newcommand{\dist}{\mathrm{dist}}

\newcommand{\GL}{\mathrm{GL}}

\newcommand{\lambdahat}{\hat{\lambda}}
\newcommand{\Mat}{\mathrm{Mat}}

\newcommand{\omegahat}{\hat{\omega}}

\newcommand{\tr}{\mathrm{tr}}

\newcommand{\uu}{\mathbf{u}}
\newcommand{\uw}{\mathbf{w}}
\newcommand{\ux}{\mathbf{x}}
\newcommand{\uy}{\mathbf{y}}
\newcommand{\uz}{\mathbf{z}}

\begin{document}

\baselineskip=17pt

\title[On two exponents of approximation]
{On two exponents of approximation related to a real number and
its square}
\author{Damien ROY}
\address{
   D\'epartement de Math\'ematiques\\
   Universit\'e d'Ottawa\\
   585 King Edward\\
   Ottawa, Ontario K1N 6N5, Canada}
\email{droy@uottawa.ca}
\thanks{Work partially supported by NSERC and CICMA}

\begin{abstract}
For each real number $\xi$, let $\lambdahat_2(\xi)$ denote the
supremum of all real numbers $\lambda$ such that, for each
sufficiently large $X$, the inequalities $|x_0| \le X$,
$|x_0\xi-x_1|  \le X^{-\lambda}$ and $|x_0\xi^2-x_2| \le
X^{-\lambda}$ admit a solution in integers $x_0$, $x_1$ and $x_2$
not all zero, and let $\omegahat_2(\xi)$ denote the supremum of
all real numbers $\omega$ such that, for each sufficiently large
$X$, the dual inequalities $|x_0+x_1\xi+x_2\xi^2| \le
X^{-\omega}$, $|x_1| \le X$ and $|x_2| \le X$ admit a solution in
integers $x_0$, $x_1$ and $x_2$ not all zero.  Answering a
question of Y.~Bugeaud and M.~Laurent, we show that the exponents
$\lambdahat_2(\xi)$ where $\xi$ ranges through all real numbers
with $[\bQ(\xi):\bQ]>2$ form a dense subset of the interval $[1/2,
(\sqrt{5}-1)/2]$ while, for the same values of $\xi$, the dual
exponents $\omegahat_2(\xi)$ form a dense subset of $[2,
(\sqrt{5}+3)/2]$. Part of the proof rests on a result of
V.~Jarn\'{\i}k showing that $\lambdahat_2(\xi) =
1-\omegahat_2(\xi)^{-1}$ for any real number $\xi$ with
$[\bQ(\xi):\bQ]>2$.
\end{abstract}

\maketitle

\section{Introduction}
 \label{intro}

Let $\xi$ and $\eta$ be real numbers.   Following the notation of
Y.~Bugeaud and M.~Laurent in \cite{BL}, we define
$\lambdahat(\xi,\eta)$ to be the supremum of all real numbers
$\lambda$ such that the inequalities
\begin{equation*}
 |x_0|\le X,
 \quad
 |x_0\xi-x_1|\le X^{-\lambda}
 \et
 |x_0\eta-x_2| \le X^{-\lambda}
\end{equation*}
admit a non-zero integer solution $(x_0, x_1, x_2)\in \bZ^3$ for
each sufficiently large value of $X$.  Similarly, we define
$\omegahat(\xi,\eta)$ to be the supremum of all real numbers
$\omega$ such that the inequalities
\begin{equation*}
 |x_0+x_1\xi+x_2\eta|\le X^{-\omega},
 \quad
 |x_1| \le X
 \et
 |x_2| \le X
\end{equation*}
admit a non-zero solution $(x_0, x_1, x_2)\in \bZ^3$ for each
sufficiently large value of $X$.  An application of Dirichlet box
principle shows that we have $1/2\le \lambdahat(\xi,\eta)$ and
$2\le \omegahat(\xi,\eta)$.  Moreover, in the (non-degenerate)
case where $1$, $\xi$ and $\eta$ are linearly independent over
$\bQ$, a result of V.~Jarn\'{\i}k, kindly pointed out to the
author by Yann Bugeaud, shows that these exponents are related by
the formula
\begin{equation}
 \label{formula:Jarnik}
 \lambdahat(\xi,\eta) = 1 - \frac{1}{\omegahat(\xi,\eta)},
\end{equation}
with the convention that the right hand side of this equality is
$1$ if $\omegahat(\xi,\eta)=\infty$ (see Theorem 1 of \cite{Ja}).

In the case where $\eta=\xi^2$, we use the shorter notation
$\lambdahat_2(\xi) := \lambdahat(\xi,\xi^2)$ and $\omegahat_2(\xi)
:= \omegahat(\xi,\xi^2)$ of \cite{BL}.  The condition that $1$,
$\xi$ and $\xi^2$ are linearly independent over $\bQ$ simply means
that $\xi$ is not an algebraic number of degree at most $2$ over
$\bQ$, a condition which we also write as $[\bQ(\xi):\bQ]>2$.
Under this condition, it is known that these exponents satisfy
\begin{equation}
 \label{ineq:exponents}
 \frac{1}{2}
 \le \lambdahat_2(\xi)
 \le \frac{1}{\gamma} = 0.618\dots
 \et
 2
 \le \omegahat_2(\xi)
 \le \gamma^2 = 2.618\dots,
\end{equation}
where $\gamma=(1+\sqrt{5})/2$ denotes the golden ratio.    By
virtue of W.~M.~Schmidt's subspace theorem, the lower bounds in
\eqref{ineq:exponents} are achieved by any algebraic number $\xi$
of degree at least 3 (see Corollaries 1C and 1E in Chapter VI of
\cite{Sc}). They are also achieved by almost all real numbers
$\xi$, with respect to Lebesgue's measure (see Theorem 2.3 of
\cite{BL}). On the other hand, the upper bounds follow
respectively from Theorem 1a of \cite{DSb} and from \cite{AR}.
They are achieved in particular by the so-called Fibonacci
continued fractions (see \S2 of \cite{Ra} or \S6 of \cite{Rb}), a
special case of the Sturmian continued fractions of \cite{ADQZ}.
Now, thanks to Jarn\'{\i}k's formula \eqref{formula:Jarnik}, we
recognize that each set of inequalities in \eqref{ineq:exponents}
can be deduced from the other one.

Generalizing the approach of \cite{Ra}, Bugeaud and Laurent have
computed the exponents $\lambdahat_2(\xi)$ and $\omegahat_2(\xi)$
for a general (characteristic) Sturmian continued fraction $\xi$.
They found that, aside from $1/\gamma$ and $\gamma^2$, the next
largest values of $\lambdahat_2(\xi)$ and $\omegahat_2(\xi)$ for
such numbers $\xi$ are respectively $2-\sqrt{2}\simeq 0.586$ and
$1+\sqrt{2}\simeq 2.414$, and they asked if there exists any
transcendental real number $\xi$ which satisfies either $2-\sqrt{2}
< \lambdahat_2(\xi) <1/\gamma$ or $1+\sqrt{2} <\omegahat_2(\xi)
<\gamma^2$ (see \S8 of \cite{BL}). Our main result below shows that
such numbers exist.

\begin{maintheorem}
The points $(\lambdahat_2(\xi),\omegahat_2(\xi))$ where $\xi$ runs
through all real numbers with $[\bQ(\xi):\bQ]>2$ form a dense
subset of the curve $\cC = \{ (1-\omega^{-1},\omega) \,;\, 2\le
\omega\le \gamma^2 \}$.
\end{maintheorem}

Since $(\lambdahat_2(\xi),\omegahat_2(\xi))=(1/2,2)$ for any
algebraic number $\xi$ of degree at least $3$, it follows in
particular that $(1/\gamma,\gamma^2)$ is an accumulation point for
the set of  points $(\lambdahat_2(\xi), \omegahat_2(\xi))$ with
$\xi$ a transcendental real number.  Because of Jarn\'{\i}k's
formula \eqref{formula:Jarnik}, this theorem is equivalent to
either one of the following two assertions.

\begin{maincorollary}
The exponents $\lambdahat_2(\xi)$ attached to transcendental real
numbers $\xi$ form a dense subset of the interval
$[1/2,1/\gamma]$. The corresponding dual exponents
$\omegahat_2(\xi)$ form a dense subset of $[2,\gamma^2]$.
\end{maincorollary}

The proof is inspired by the constructions of \S6 of \cite{Rb} and
\S5 of \cite{Rd}.  We produce countably many real numbers $\xi$ of
``Fibonacci type'' (see \S\ref{section:exponents} for a precise
definition) for which we show that the exponents
$\omegahat_2(\xi)$ are dense in $[2,\gamma^2]$.  By
\eqref{formula:Jarnik}, this implies the theorem. One may then
reformulate the question of Bugeaud and Laurent by asking if there
exist transcendental real numbers $\xi$ not of that type which
satisfy $\omegahat_2(\xi) > 1+\sqrt{2}$. The work of S.~Fischler
announced in \cite{Fi} should bring some light on this question.

\begin{acknowlegments}
The author warmly thanks Yann Bugeaud for pointing out the results
of Jarn\'{\i}k in \cite{Ja} which brought a notable simplification
to the present paper.
\end{acknowlegments}

\section{Notation and equivalent definitions of the exponents}
 \label{section:notation}

We define the {\it norm} of a point $\ux=(x_0,x_1,x_2)\in\bR^3$ as
its maximum norm
\begin{equation*}
 \|\ux\|=\max_{0\le i\le 2} |x_i|.
\end{equation*}
Given a second point $\uy\in\bR^3$, we denote by $\ux\wedge \uy$
the standard vector product of $\ux$ and $\uy$, and by $\langle
\ux, \uy \rangle$ their standard scalar product. Given a third
point $\uz\in\bR^3$, we also denote by $\det(\ux,\uy,\uz)$ the
determinant of the $3\times 3$ matrix whose rows are $\ux$, $\uy$
and $\uz$.  Then we have the well-known relation
\begin{equation*}
 \det(\ux,\uy,\uz)=\langle \ux, \uy\wedge \uz\rangle
\end{equation*}
and we get the following alternative definition
of the exponents $\lambdahat(\xi,\eta)$ and $\omegahat(\xi,\eta)$.

\begin{lemma}
 \label{lemma:exponents}
Let $\xi,\eta\in\bR$, and let $\uy=(1,\xi,\eta)$.  Then
$\lambdahat(\xi,\eta)$ is the supremum of all real numbers
$\lambda$ such that, for each sufficiently large real number $X\ge
1$, there exists a point $\ux\in\bZ^3$ with
\begin{equation*}
 0<\|\ux\|\le X
 \et
 \|\ux\wedge \uy\| \le X^{-\lambda}.
\end{equation*}
Similarly, $\omegahat(\xi,\eta)$ is the supremum of all real
numbers $\omega$ such that, for each sufficiently large real
number $X\ge 1$, there exists a point $\ux\in\bZ^3$ with
\begin{equation*}
 0<\|\ux\|\le X
 \et
 |\langle \ux, \uy \rangle| \le X^{-\omega}.
\end{equation*}
\end{lemma}

\medskip
In the sequel, we will need the following inequalities.

\begin{lemma}
 \label{lemma:ineq}
For any $\ux,\uy,\uz\in\bR^3$, we have
\begin{equation}
 \label{ineq:pscal}
 \|\langle \ux, \uz \rangle \uy - \langle \ux, \uy \rangle \uz \|
  \le 2\|\ux\| \|\uy \wedge \uz\|,
\end{equation}
\begin{equation}
 \label{ineq:pvect}
 \|\uy\| \| \ux \wedge \uz \|
  \le \|\uz\| \| \ux \wedge \uy \| + 2\|\ux\| \|\uy \wedge \uz\|.
\end{equation}
\end{lemma}

\begin{proof} Writing $\uy=(y_0,y_1,y_2)$ and $\uz=(z_0,z_1,z_2)$, we
find
\begin{equation*}
 \| \langle \ux, \uz \rangle \uy - \langle \ux, \uy \rangle \uz \|
  = \max_{i=0,1,2} | \langle \ux, y_i\uz-z_i\uy \rangle|
  \le 2\|\ux\| \|\uy \wedge \uz\|,
\end{equation*}
which proves \eqref{ineq:pscal}.  Similarly, one finds $\| y_i \ux
\wedge \uz - z_i \ux \wedge \uy \| \le 2\|\ux\| \|\uy \wedge \uz\|$
for $i=0,1,2$, and this implies \eqref{ineq:pvect}.
\end{proof}

For any non-zero point $\ux$ of $\bR^3$, let $[\ux]$ denote the
point of $\bP^2(\bR)$ having $\ux$ as a set of homogeneous
coordinates. Then, \eqref{ineq:pvect} has a useful interpretation
in terms of the projective distance defined for non-zero points
$\ux$ and $\uy$ of $\bR^3$ by
\begin{equation*}
 \dist([\ux],[\uy]) = \dist(\ux,\uy) = \frac{\|\ux\wedge \uy\|}{\|\ux\| \|\uy\|}.
\end{equation*}
Indeed, for any triple of non-zero points $\ux,\uy,\uz\in\bR^3$,
it gives
\begin{equation}
 \label{ineq:dist}
 \dist( [\ux], [\uz] ) \le \dist( [\ux], [\uy] ) + 2\, \dist( [\uy], [\uz] ).
\end{equation}

\section{Fibonacci sequences in $\GL_2(\bC)$}
\label{section:Fibonacci}

A {\it Fibonacci sequence} in a monoid is a sequence
$(\uw_i)_{i\ge 0}$ of elements of that monoid such that
$\uw_{i+2}=\uw_{i+1}\uw_i$ for each index $i\ge 0$.  Clearly, such
a sequence is entirely determined by its first two elements
$\uw_0$ and $\uw_1$.  We start with the following observation.

\begin{proposition}
 \label{prop:GL}
There exists a non-empty Zariski open subset $\cU$ of
$\GL_2(\bC)^2$ with the following property.  For each Fibonacci
sequence $(\uw_i)_{i\ge 0}$ with $(\uw_0,\uw_1)\in \cU$, there
exists $N\in \GL_2(\bC)$ such that the matrix
\begin{equation}
 \label{sequences:yi}
 \uy_i =
 \begin{cases}
  \uw_iN &\mbox{if $i$ is even,}\\
  \uw_i{^tN} &\mbox{if $i$ is odd}
 \end{cases}
\end{equation}
is symmetric for each $i\ge 0$.  Any matrix $N\in\GL_2(\bC)$ such
that $\uw_0N$, $\uw_1{^tN}$ and $\uw_1\uw_0N$ are symmetric
satisfies this property.  When $\uw_0$ and $\uw_1$ have integer
coefficients, we may take $N$ with integer coefficients.
\end{proposition}

\begin{proof}  Let $(\uw_i)_{i\ge 0}$ be a Fibonacci sequence in
$\GL_2(\bC)$ and let $N\in\GL_2(\bC)$.  Defining $\uy_i$ by
\eqref{sequences:yi} for each $i\ge 0$, we find $\uy_{i+3} =
\uy_{i+1} {^tS} \uy_i S \uy_{i+1}$ with $S=N^{-1}$ if $i$ is even
and $S={^t}N^{-1}$ if $i$ is odd.  Thus, $\uy_i$ is symmetric for
each $i\ge 0$ if and only if it is so for $i=0,1,2$.

Now, for any given point $(\uw_0,\uw_1)\in \GL_2(\bC)^2$, the
conditions that $\uw_0N$, $\uw_1{^tN}$ and $\uw_1\uw_0N$ are
symmetric represent a system of three linear equations in the four
unknown coefficients of $N$.  Let $\cV$ be the Zariski open subset
of $\GL_2(\bC)^2$ consisting of all points $(\uw_0,\uw_1)$ for
which this linear system has rank $3$. Then, for each
$(\uw_0,\uw_1)\in \cV$, the $3\times 3$ minors of this linear
system conveniently arranged into a $2\times 2$ matrix provide a
non-zero solution $N$ of the system, whose coefficients are
polynomials in those of $\uw_0$ and $\uw_1$ with integer
coefficients. Then, the condition $\det(N)\neq 0$ in turn
determines a Zariski open subset $\cU$ of $\cV$.  To conclude, we
note that $\cU$ is not empty as a short computation shows that it
contains the point formed by $\uw_0=\begin{pmatrix} 1 &1\\ 1
&0\end{pmatrix}$ and $\uw_1=\begin{pmatrix} 0 &1\\ 1
&0\end{pmatrix}$.
\end{proof}

\begin{definition}
 \label{def1}
Let $\cM=\Mat_{2\times 2}(\bZ)\cap \GL_2(\bC)$ denote the monoid
of $2\times 2$ integer matrices with non-zero determinant. We say
that a Fibonacci sequence $(\uw_i)_{i\ge 0}$ in $\cM$ is {\it
admissible} if there exists a matrix $N\in\cM$ such that the
sequence $(\uy_i)_{i\ge 0}$ given by \eqref{sequences:yi} consists
of symmetric matrices.
\end{definition}

Since $\cM$ is Zariski dense in $\GL_2(\bC)$, Proposition
\ref{prop:GL} shows that almost all Fibonacci sequences in
$\cM$ are admissible.  The following example is an illustration
of this.

\begin{example}
 \label{example1}
Fix integers $a,b,c$ with $a\ge 2$ and $c\ge b\ge 1$, and define
\begin{equation*}
 \uw_0 = \begin{pmatrix} 1 &b\\ a &a(b+1)\end{pmatrix},
 \quad
 \uw_1 = \begin{pmatrix} 1 &c\\ a &a(c+1)\end{pmatrix}
\end{equation*}
and
\begin{equation*}
 N = \begin{pmatrix} -1+a(b+1)(c+1) &-a(b+1) \\ -a(c+1) &a
     \end{pmatrix}.
\end{equation*}
These matrices belong to $\cM$ since $\det(\uw_0)=\det(\uw_1)=a$
and $\det(N)=-a$.  Moreover, one finds that
\begin{equation*}
 \uw_0N = \begin{pmatrix} -1+a(c+1) &-a\\ -a &0\end{pmatrix},
 \quad
 \uw_1{^t}N = \begin{pmatrix} -1+a(b+1) &-a\\ -a &0\end{pmatrix}
\end{equation*}
and
\begin{equation*}
 \uw_1\uw_0N = \begin{pmatrix}-1+a &-a\\ -a &-a^2\end{pmatrix}
\end{equation*}
are symmetric matrices. Therefore, the Fibonacci sequence
$(\uw_i)_{i\ge0}$ constructed on $\uw_0$ and $\uw_1$ is admissible
with an associated sequence of symmetric matrices $(\uy_i)_{i\ge0}$
given by \eqref{sequences:yi}, the first three matrices of this
sequence being the above products $\uy_0=\uw_0N$, $\uy_1=\uw_1{^t}N$
and $\uy_2=\uw_1\uw_0N$.
\end{example}

\section{Fibonacci sequences of $2\times 2$ integer matrices}
\label{section:arithmetic}

In the sequel, we identify $\bR^3$ (resp.\ $\bZ^3$) with the space
of $2\times 2$ symmetric matrices with real (resp.\ integer)
coefficients under the map
\begin{equation*}
 \ux=(x_0,x_1,x_2)
 \ \longmapsto\
 \begin{pmatrix} x_0 &x_1\\x_1 &x_2 \end{pmatrix}.
\end{equation*}
Accordingly, it makes sense to define the determinant of a point
$\ux=(x_0,x_1,x_2)$ of $ \bR^3$ by $\det(\ux) = x_0x_2-x_1^2$.
Similarly, given symmetric matrices $\ux$, $\uy$ and $\uz$, we write
$\ux\wedge \uy$, $\langle \ux, \uy\rangle$ and $\det(\ux,\uy,\uz)$
to denote respectively the vector product, scalar product and
determinant of the corresponding points.

In this section we look at arithmetic properties of admissible
Fibonacci sequences in the monoid $\cM$ of Definition \ref{def1}.
For this purpose, we define the {\it content} of an integer matrix
$\uw \in \Mat_{2\times 2}(\bZ)$ or of a point $\uy \in \bZ^3$ as the
greatest common divisor of their coefficients. We say that such a
matrix or point is {\it primitive} if its content is $1$.

\begin{proposition}
 \label{prop:formulesrec}
Let $(\uw_i)_{i\ge 0}$ be an admissible Fibonacci sequence of
matrices in $\cM$ and let $(\uy_i)_{i\ge 0}$ be a corresponding
sequence of symmetric matrices in $\cM$.  For each $i\ge 0$,
define $\uz_i=\det(\uw_i)^{-1}\uy_i\wedge \uy_{i+1}$. Then, for
each $i\ge 0$, we have
\begin{itemize}
\item[(a)]
 $\tr(\uw_{i+3})=\tr(\uw_{i+1})\tr(\uw_{i+2})-\det(\uw_{i+1})\tr(\uw_i)$,
\item[(b)]
 $\uy_{i+3}=\tr(\uw_{i+1})\uy_{i+2}-\det(\uw_{i+1})\uy_i$,
\item[(c)]
 $\uz_{i+3}=\tr(\uw_{i+1})\uz_{i+1}+\det(\uw_i)\uz_i$,
\item[(d)]
 $\det(\uy_i,\uy_{i+1},\uy_{i+2}) =
 (-1)^i\det(\uy_0,\uy_1,\uy_2)\det(\uw_2)^{-1}\det(\uw_{i+2})$,
\item[(e)]
 $\uz_i\wedge \uz_{i+1} =
 (-1)^i \det(\uy_0,\uy_1,\uy_2)\det(\uw_2)^{-1}\uy_{i+1}$.
\end{itemize}
\end{proposition}

\begin{proof}
For each index $i\ge 0$, let $N_i$ denote the element of $\cM$ for
which $\uy_i=\uw_iN_i$.  According to \eqref{sequences:yi}, we have
$N_i=N$ if $i$ is even and $N_i={^t}N$ if $i$ is odd.  We first
prove (b) following the argument of the proof of Lemma 2.5 (i) of
\cite{Rc}.  Multiplying both sides of the equality
$\uw_{i+2}=\uw_{i+1}\uw_i$ on the right by $N_{i+2}=N_i$, we find
\begin{equation}
 \label{facto:wy}
 \uy_{i+2}=\uw_{i+1}\uy_i,
\end{equation}
which can be rewritten as $\uy_{i+2}=\uy_{i+1}N_{i+1}^{-1}\uy_i$.
Taking transpose of both sides, this gives $\uy_{i+2}=\uy_i
N_i^{-1}\uy_{i+1}=\uw_i\uy_{i+1}$. Replacing $i$ by $i+1$ in the
latter identity and combining it with \eqref{facto:wy}, we get
\begin{equation}
 \label{facto:wwy}
 \uy_{i+3}=\uw_{i+1}\uy_{i+2}=\uw_{i+1}^2\uy_i.
\end{equation}
Then (b) follows from \eqref{facto:wy} and \eqref{facto:wwy} using
the fact that, by the Cayley-Hamilton theorem, we have $\uw_{i+1}^2
= \tr(\uw_{i+1})\uw_{i+1} - \det(\uw_{i+1})I$.  Multiplying both
sides of (b) on the right by $N_i^{-1}$ and taking the trace, we
deduce that
\begin{equation*}
 \tr(\uy_{i+3}N_i^{-1})
 = \tr(\uw_{i+1})\tr(\uw_{i+2})-\det(\uw_{i+1})\tr(\uw_i).
\end{equation*}
This gives (a) because $\tr(\uy_{i+3}N_i^{-1}) =
\tr({^t}\uy_{i+3}{^t}N_i^{-1}) = \tr(\uw_{i+3})$.  Taking the
exterior product of both sides of (b) with $\uy_{i+1}$, we also find
\begin{equation*}
 \uy_{i+1}\wedge \uy_{i+3}
  =
 \tr(\uw_{i+1})\det(\uw_{i+1})\uz_{i+1}
       +\det(\uw_{i+1})\det(\uw_i)\uz_i.
\end{equation*}
Similarly, replacing $i$ by $i+1$ in (b) and taking the exterior
product with $\uy_{i+3}$ gives
\begin{equation*}
 \det(\uw_{i+3})\uz_{i+3}
  = \det(\uw_{i+2}) \uy_{i+1}\wedge \uy_{i+3}.
\end{equation*}
Then (c) follows upon noting that $\det(\uw_{i+3}) =
\det(\uw_{i+2}) \det(\uw_{i+1})$.

The formula (d) is clearly true for $i=0$.  If we assume that it
holds for some integer $i\ge 0$, then using the formula for
$\uy_{i+3}$ given by (b) and taking into account the
multilinearity of the determinant we find
\begin{equation*}
 \det(\uy_{i+1},\uy_{i+2},\uy_{i+3})
 = -\det(\uw_{i+1})\det(\uy_i,\uy_{i+1},\uy_{i+2})
 = (-1)^{i+1}\det(\uy_0,\uy_1,\uy_2)\frac{\det(\uw_{i+3})}{\det(\uw_2)}.
\end{equation*}
This proves (d) by induction on $i$.  Then (e) follows since, for
any $\ux,\uy,\uz\in\bZ^3$, we have $(\ux\wedge
\uy)\wedge(\uy\wedge \uz)=\det(\ux,\uy,\uz)\,\uy$ which, in the
present case, gives
\begin{equation*}
 \uz_i\wedge \uz_{i+1} =
 \det(\uw_{i+2})^{-1}\det(\uy_i,\uy_{i+1},\uy_{i+2})\,\uy_{i+1}.
\end{equation*}
\end{proof}

\begin{corollary}
 \label{corollary:arith}
The notation being as in the proposition, assume that $\tr(\uw_i)$
and $\det(\uw_i)$ are relatively prime for $i=0,1,2,3$ and that
$\det(\uy_0,\uy_1,\uy_2)\neq 0$.  Then, for each $i\ge 0$,
\begin{itemize}
 \item[(a)] the points $\uy_i,\uy_{i+1},\uy_{i+2}$ are linearly
 independent,
 \item[(b)] $\tr(\uw_i)$ and $\det(\uw_i)$ are relatively prime,
 \item[(c)] the matrix $\uw_i$ is primitive,
 \item[(d)] the content of $\uy_i$ divides $\det(\uy_2)/\det(\uw_2)$,
 \item[(e)] the point $\det(\uw_2)\,\uz_i$ belongs to $\bZ^3$ and its
 content divides $\det(\uy_2)\det(\uy_0,\uy_1,\uy_2)$.
\end{itemize}
\end{corollary}

\begin{proof}
The assertion (a) follows from Proposition \ref{prop:formulesrec}
(d). Since (b) holds by hypothesis for $i=0,1,2,3$, and since
$\det(\uw_2)$ and $\det(\uw_i)$ have the same prime factors for
each $i\ge 2$, the assertion (b) follows, by induction on $i$,
from the fact that Proposition \ref{prop:formulesrec} (a) gives
$\tr(\uw_{i+1})\equiv \tr(\uw_i)\tr(\uw_{i-1})$ modulo
$\det(\uw_2)$ for each $i\ge 3$. Then (c) follows since the
content of $\uw_i$ divides both $\tr(\uw_i)$ and $\det(\uw_i)$.

Let $N\in\cM$ such that $\uy_2=\uw_2N$.  For each $i$, we have
$\uy_i=\uw_iN_i$ where $N_i=N$ if $i$ is even and $N_i={^t}N$ if $i$
is odd.  This gives $\uy_i\adj(N_i) = \det(N)\uw_i$ where
$\adj(N_i)\in\cM$ denotes the adjoint of $N_i$.  Thus, by (c), the
content of $\uy_i$ divides $\det(N) = \det(\uy_2)/\det(\uw_2)$, as
claimed in (d).

The fact that $\det(\uw_2)\,\uz_i$ belongs to $\bZ^3$ is clear for
$i=0,1,2$ because $\det(\uw_0)$ and $\det(\uw_1)$ divide
$\det(\uw_2)$. Then, Proposition \ref{prop:formulesrec} (c) shows,
by induction on $i$, that $\det(\uw_2)\,\uz_i\in\bZ^3$ for each
$i\ge 0$. Moreover, the content of that point divides that of
$\det(\uw_2)^2\uz_i\wedge \uz_{i+1}$ which, by (d) and Proposition
\ref{prop:formulesrec} (e), divides
$\det(\uy_0,\uy_1,\uy_2)\det(\uy_2)$. This proves (e).
\end{proof}

\begin{example}
 \label{example2}
Let $(\uw_i)_{i\ge 0}$, $N$ and $(\uy_i)_{i\ge 0}$ be as in
Example \ref{example1}.  Since $\uw_0$, $\uw_1$ and $N$ are
congruent to matrices of the form $\begin{pmatrix} \pm 1 &*\\
0&0\end{pmatrix}$ modulo $a$ and have determinant $\pm a$, all
matrices $\uw_i$ and $\uy_i$ are congruent to matrices of the same
form modulo $a$ and their determinant is, up to sign, a power of
$a$. Thus these matrices have relatively prime trace and
determinant, and so are primitive for each $i\ge 0$. Since
$\det(\uy_0,\uy_1,\uy_2)=a^4(c-b)$, Proposition
\ref{prop:formulesrec} (e) shows that the points
$\uz_i=\det(\uw_i)^{-1}\uy_i\wedge \uy_{i+1}$ satisfy $\uz_i\wedge
\uz_{i+1}=(-1)^i a^2 (c-b) \uy_{i+1}$ for each $i\ge 0$. Moreover,
we find that $a^{-1}\uz_0=(0,0,b-c)$,
$a^{-1}\uz_1=(a,-1+a(b+1),-b)$ and $a^{-1}\uz_2=(a,-1+a(c+1),-c)$
are integer points.  Then, Proposition \ref{prop:formulesrec} (c)
shows, by induction on $i$, that $a^{-1}\uz_i\in\bZ^3$ for each
$i\ge 0$. In particular, if $c=b+1$, we deduce from the relation
$a^{-1}\uz_i\wedge a^{-1}\uz_{i+1}= \pm \uy_{i+1}$ that
$a^{-1}\uz_i$ is a primitive integer point for each $i\ge 0$.
\end{example}

\section{Growth estimates}
\label{section:growth}

Define the {\it norm} of a $2\times 2$ matrix $\uw =( w_{k,\ell})
\in \Mat_{2\times 2}(\bR)$ as the largest absolute value of its
coefficients $\|\uw\|=\max_{1\le k,\ell\le 2}|w_{k,\ell}|$, and
define $\gamma=(1+\sqrt{5})/2$ as in the introduction. In this
section, we provide growth estimates for the norm and determinant of
elements of certain Fibonacci sequences in $\GL_2(\bR)$.  We first
establish two basic lemmas.

\begin{lemma}
 \label{estimates:lemma:special}
Let $\uw_0,\uw_1\in\GL_2(\bR)$.  Suppose that, for $i=0,1$, the
matrix $\uw_i$ is of the form
\begin{equation*}
 \begin{pmatrix} a&b\\c&d \end{pmatrix}
 \mbox{\quad with \quad}
 1 \le a \le \min\{b,c\}
 \et
 \max\{b,c\} \le d.
\end{equation*}
Then, all matrices of the Fibonacci sequence $(\uw_i)_{i\ge 0}$
constructed on $\uw_0$ and $\uw_1$ have this form and, for each
$i\ge 0$, they satisfy
\begin{equation}
 \label{estimates:growth_pos}
 \|\uw_i\|\|\uw_{i+1}\| < \|\uw_{i+2}\| \le 2\|\uw_i\|\|\uw_{i+1}\|.
\end{equation}
\end{lemma}

\begin{proof}
The first assertion follows by recurrence on $i$ and is left to the
reader. It implies that $\|\uw_i\|$ is equal to the element of index
$(2,2)$ of $\uw_i$ for each $i\ge 0$. Then,
\eqref{estimates:growth_pos} follows by observing that, for any
$2\times 2$ matrices $\uw=(w_{k,\ell})$ and $\uw'=(w'_{k,\ell})$
with positive real coefficients, the product
$\uw'\uw=(w''_{k,\ell})$ satisfies $w_{2,2}w'_{2,2} < w''_{2,2} \le
2\|\uw\|\|\uw'\|$.
\end{proof}

\begin{lemma}
 \label{estimates:lemma:numbers}
Let $(r_i)_{i\ge 0}$ be a sequence of positive real numbers.
Assume that there exist constants $c_1,c_2>0$ such that
$c_1r_ir_{i+1} \le r_{i+2} \le c_2r_ir_{i+1}$ for each $i\ge 0$.
Then there also exist constants $c_3,c_4>0$ such that
$c_3r_i^\gamma \le r_{i+1} \le c_4r_i^\gamma$ for each $i\ge 0$.
\end{lemma}

\begin{proof}
Define $c_3=c_1^\gamma/(cc_2)$ and $c_4=cc_2^\gamma/c_1$ where
$c\ge 1$ is chosen so that the condition $c_3 \le
r_{i+1}/r_i^\gamma \le c_4$ holds for $i=0$.  Assuming that the
same condition holds for some index $i\ge 0$, we find
\begin{equation*}
 \frac{r_{i+2}}{r_{i+1}^\gamma}
  \ge c_1 \frac{r_i}{r_{i+1}^{1/\gamma}}
  \ge c_1c_4^{-1/\gamma}
  = c^{1/\gamma^2}c_3
  \ge c_3,
\end{equation*}
and similarly $r_{i+2}/r_{i+1}^\gamma \le c_4$.  This proves the
lemma by recurrence on $i$.
\end{proof}

\begin{proposition}
 \label{prop:growth}
Let $(\uw_i)_{i\ge 0}$ be a Fibonacci sequence in $\GL_2(\bR)$.
Suppose that there exist real numbers $c_1,c_2>0$ such that
\begin{equation}
 \label{estimates:growth}
 c_1 \|\uw_i\|\|\uw_{i+1}\| \le \|\uw_{i+2}\| \le c_2\|\uw_i\|\|\uw_{i+1}\|
\end{equation}
for each $i\ge 0$.  Then, there exist constants $c_3,c_4>0$ such
that
\begin{equation}
 \label{estimates:norm,det}
 c_3\|\uw_i\|^\gamma \le \|\uw_{i+1}\| \le c_4\|\uw_i\|^\gamma
 \et
 c_3|\det(\uw_i)|^\gamma \le |\det(\uw_{i+1})| \le c_4|\det(\uw_i)|^\gamma
\end{equation}
for each $i\ge 0$.  Moreover, if there exist $\alpha,\beta\ge 0$
such that
\begin{equation}
 \label{estimates:det}
 (c_2\|\uw_i\|)^\alpha \le |\det(\uw_i)| \le (c_1\|\uw_i\|)^\beta
\end{equation}
holds for $i=0,1$, then the latter relation extends to each $i\ge
0$.
\end{proposition}

\begin{proof}
The first assertion of the proposition follows from Lemma
\ref{estimates:lemma:numbers} applied once with $r_i=\|\uw_i\|$ and
once with $r_i=|\det(\uw_i)|$.  To prove the second assertion,
assume that, for some index $j\ge 0$, the condition
\eqref{estimates:det} holds both with $i=j$ and $i=j+1$. We find
\begin{equation*}
 |\det(\uw_{j+2})|
  =|\det(\uw_{j+1})||\det(\uw_{j})|
  \ge (c_2\|\uw_{j+1}\|)^\alpha (c_2\|\uw_{j}\|)^\alpha
  \ge (c_2\|\uw_{j+2}\|)^\alpha
\end{equation*}
and similarly $|\det(\uw_{j+2})| \le (c_1\|\uw_{j+2}\|)^\beta$.
Therefore, \eqref{estimates:det} holds with $i=j+2$.  By recurrence
on $i$, this shows that \eqref{estimates:det} holds for each $i\ge
0$ if it holds for $i=0,1$.
\end{proof}

\begin{example}
 \label{example3}
Let the notation be as in Example \ref{example1}.  Since $\uw_0$ and
$\uw_1$ satisfy the hypotheses of Lemma
\ref{estimates:lemma:special}, the Fibonacci sequence $(\uw_i)_{i\ge
0}$ that they generate fulfills for each $i\ge 0$ the condition
\eqref{estimates:growth} of Proposition \ref{prop:growth} with
$c_1=1$ and $c_2=2$.  As $\det(\uw_0)=\det(\uw_1)=a$, we also note
that, for this choice of $c_1$ and $c_2$, the condition
\eqref{estimates:det} holds for $i=0,1$ with
\begin{equation*}
 \alpha=\frac{\log a}{\log(2a(c+1))}
 \et
 \beta=\frac{\log a}{\log(a(b+1))}.
\end{equation*}
Then, for an appropriate choice of $c_3,c_4>0$, both
\eqref{estimates:norm,det} and \eqref{estimates:det} hold for each
$i\ge 0$.  Moreover, the estimates \eqref{estimates:growth_pos} of
Lemma \ref{estimates:lemma:special} imply that the sequence
$(\uw_i)_{i\ge 0}$ is unbounded.
\end{example}

\section{Construction of a real number}
\label{section:construction}

Given sequences of non-negative real numbers with general terms
$a_i$ and $b_i$, we write $a_i\ll b_i$ or $b_i \gg a_i$ if there
exists a real number $c>0$ such that $a_i\le cb_i$ for all
sufficiently large values of $i$.  We write $a_i\sim b_i$ when
$a_i\ll b_i$ and $b_i\ll a_i$.  With this notation, we now prove
the following result (compare with \S5 of \cite{Rd}).

\begin{proposition}
 \label{prop:convergence}
Let $(\uw_i)_{i\ge0}$ be an admissible Fibonacci sequence in $\cM$
and let $(\uy_i)_{i\ge 0}$ be a corresponding sequence of symmetric
matrices in $\cM$.  Assume that $(\uw_i)_{i\ge0}$ is unbounded and
satisfies the conditions
\begin{equation}
 \label{ineq:conv0}
 \|\uw_{i+1}\| \sim \|\uw_i\|^\gamma,
 \quad
 |\det(\uw_{i+1})| \sim |\det(\uw_i)|^\gamma
 \et
 |\det(\uw_i)| \ll \|\uw_i\|^\beta
\end{equation}
for a real number $\beta$ with $0<\beta<2$.  Viewing each $\uy_i$ as
a point in $\bZ^3$, assume that $\det(\uy_0,\uy_1,\uy_2)\neq 0$ and
define $\uz_i=(\det(\uw_i))^{-1}\uy_i\wedge \uy_{i+1}$ for each
$i\ge 0$. Then we have
\begin{equation}
 \label{ineq:conv1}
 \|\uy_i\| \sim \|\uw_i\|,
 \quad
 |\det(\uy_i)| \sim|\det(\uw_i)|,
 \quad
 \|\uz_i\| \sim \|\uw_{i-1}\|,
\end{equation}
and there exists a non-zero point $\uy$ of $\bR^3$ with
$\det(\uy)=0$ such that
\begin{equation}
 \label{ineq:conv2}
 \|\uy_i\wedge \uy\| \sim \frac{|\det(\uw_i)|}{\|\uw_i\|}
 \et
 |\langle \uz_i,\uy \rangle| \sim
 \frac{|\det(\uw_{i+1})|}{\|\uw_{i+2}\|}.
\end{equation}
If $\beta<1$, the coordinates of such a point $\uy$ are linearly
independent over $\bQ$ and we may assume that $\uy=(1,\xi,\xi^2)$
for some real number $\xi$ with $[\bQ(\xi):\bQ]
> 2$.
\end{proposition}

\begin{proof}
For each $i\ge 0$, let $N_i$ denote the element of $\cM$ for which
$\uy_i=\uw_iN_i$. Putting $N=N_0$, we have by hypothesis $N_i=N$
when $i$ is even and $N_i={^t}N$ otherwise.  This implies that
$\|\uy_i\|\sim \|\uw_i\|$ and $|\det(\uy_i)|\sim|\det(\uw_i)|$. In
the sequel, we will repeatedly use these relations as well as the
hypothesis \eqref{ineq:conv0}.

We claim that we have
\begin{equation}
 \label{ineq:wedgeyy}
 \|\uy_i\wedge \uy_{i+1}\| \ll |\det(\uw_i)| \|\uw_{i-1}\|.
\end{equation}
To prove this, we define $J = \begin{pmatrix} 0 &1\\ -1 &0
\end{pmatrix}$ and note that, for each $i\ge 0$, the coefficients
of the diagonal of $\uy_iJ\uy_{i+1}$ coincide with the first and
third coefficients of $\uy_i\wedge \uy_{i+1}$ while the sum of the
coefficients of $\uy_iJ\uy_{i+1}$ outside of the diagonal is the
middle coefficient of $\uy_i\wedge \uy_{i+1}$ multiplied by $-1$.
This gives
\begin{equation}
 \label{ineq:J}
 \|\uy_i\wedge \uy_{i+1}\| \le 2\|\uy_iJ \uy_{i+1}\|.
\end{equation}
Since $\uy_{i+1}=\uy_iN_i^{-1}\uy_{i-1}$ and since $\ux J \ux =
\det(\ux)J$ for any symmetric matrix $\ux$, we also find that
$\uy_i J \uy_{i+1} = \det(\uy_i)JN_i^{-1}\uy_{i-1}$ and therefore
$\|\uy_i J \uy_{i+1}\| \ll 
|\det(\uw_i)|\|\uw_{i-1}\|$. Combining this with \eqref{ineq:J}
proves our claim \eqref{ineq:wedgeyy} which can also be written in
the form
\begin{equation}
 \label{upperbound:normz}
 \|\uz_i\| \ll \| \uw_{i-1}\|.
\end{equation}
As $\|\uy_i\|\sim \|\uw_i\|$ and $\|\uy_{i+1}\|\sim
\|\uw_i\|^\gamma$, the estimate \eqref{ineq:wedgeyy} shows, in the
notation of \S\ref{section:notation}, that
\begin{equation}
 \label{upperbound:distyy}
 \dist([\uy_i],[\uy_{i+1}]) \le c \delta_i,
 \quad
 \mbox{where}
 \quad
 \delta_i=\frac{|\det(\uw_i)|}{\|\uw_i\|^2}
\end{equation}
and where $c$ is some positive constant which does not depend on
$i$. Since by hypothesis we have $|\det(\uw_i)|\ll
\|\uw_i\|^\beta$ with $\beta<2$, we find that
$\lim_{i\to\infty}\delta_i=0$.  Since moreover, we have
$\delta_{i+1}\sim\delta_i^\gamma$, we deduce that there exists an
index $i_0\ge 1$ such that $\delta_{i+1}\le \delta_i/4$ for each
$i\ge i_0$.  Then, using \eqref{ineq:dist}, we deduce that
\begin{equation}
 \label{ineq:distyy}
 \dist([\uy_i],[\uy_j])
 \le
 \sum_{k=i}^{j-1} 2^{k-i}\dist([\uy_k],[\uy_{k+1}])
 \le
 c \sum_{k=i}^{j-1} 2^{k-i}\delta_k
 \le
 2c\delta_i
\end{equation}
for each choice of $i$ and $j$ with $i_0\le i<j$.  Thus the
sequence $([\uy_i])_{i\ge 0}$ converges in $\bP^2(\bR)$ to a point
$[\uy]$ for some non-zero $\uy\in\bR^3$.  Since the ratio
$|\det(\uy_i)|/\|\uy_i\|^2$ depends only on the class $[\uy_i]$ of
$\uy_i$ in $\bP^2(\bR)$ and tends to $0$ like $\delta_i$ as $i\to
\infty$, we deduce, by continuity, that $|\det(\uy)|/\|\uy\|^2=0$
and thus that $\det(\uy)=0$. By continuity, \eqref{ineq:distyy}
also leads to $\dist([\uy_i],[\uy]) \le 2c\delta_i$ for each $i\ge
i_0$, and so
\begin{equation}
 \label{upperbound:wedgey}
 \|\uy_i\wedge \uy\| \ll \frac{|\det(\uw_i)|}{\|\uw_i\|}.
\end{equation}
Applying \eqref{ineq:pscal} together with the above estimates
\eqref{upperbound:normz} and \eqref{upperbound:wedgey}, we find
\begin{equation*}
 \| \langle \uz_i, \uy \rangle \uy_{i+2}
       - \langle \uz_i, \uy_{i+2} \rangle \uy \|
   \le 2\|\uz_i\| \|\uy_{i+2} \wedge \uy\|
   \ll \|\uw_{i-1}\|\frac{|\det(\uw_{i+2})|}{\|\uw_{i+2}\|}
   \ll |\det(\uw_{i+1})|\delta_i.
\end{equation*}
Using Proposition \ref{prop:formulesrec} (d), we also get
\begin{equation}
 \label{ineq:det}
 \|\langle \uz_i, \uy_{i+2} \rangle \uy \|
 =  \frac{|\det(\uy_i,\uy_{i+1},\uy_{i+2})|}{|\det(\uw_i)|} \|\uy\|
 \sim |\det(\uw_{i+1})|.
\end{equation}
Combining the above two estimates, we deduce that $\|\langle
\uz_i, \uy \rangle \uy_{i+2}\| \sim |\det(\uw_{i+1})|$ and
therefore that $|\langle \uz_i, \uy \rangle| \sim
|\det(\uw_{i+1})|/\|\uw_{i+2}\|$. The latter estimate is the
second half of \eqref{ineq:conv2}.  It implies
\begin{equation*}
 \| \langle \uz_{i+1}, \uy \rangle \uy_i \|
  \sim \frac{|\det(\uw_{i+2})|}{\|\uw_{i+3}\|} \|\uw_i\|
  \sim |\det(\uw_i)| \delta_{i+1}.
\end{equation*}
Since $\langle \uz_{i+1}, \uy_i \rangle = \det(\uw_{i-1})^{-1}
\langle \uz_i, \uy_{i+2} \rangle$, the estimate \eqref{ineq:det}
can also be written in the form $\| \langle \uz_{i+1}, \uy_i
\rangle \uy \| \sim |\det(\uw_i)|$. Then, applying
\eqref{ineq:pscal} once again, we find
\begin{equation*}
 2\|\uz_{i+1}\| \|\uy_i \wedge \uy\|
   \ge \| \langle \uz_{i+1}, \uy \rangle \uy_i
          - \langle \uz_{i+1}, \uy_i \rangle \uy \|
   \gg |\det(\uw_i)|.
\end{equation*}
Since, by \eqref{upperbound:normz} and \eqref{upperbound:wedgey},
we have $\|\uz_{i+1}\|\ll \|\uw_i\|$ and $\|\uy_i \wedge \uy\|\ll
|\det(\uw_i)|/\|\uw_i\|$, we conclude from this that
$\|\uz_{i+1}\|\sim \|\uw_i\|$ and $\|\uy_i \wedge \uy\|\sim
|\det(\uw_i)|/\|\uw_i\|$, which completes the proof of
\eqref{ineq:conv1} and \eqref{ineq:conv2}.

Now, assume that $\beta<1$, and let $\uu\in\bZ^3$ such that
$\langle \uu, \uy \rangle = 0$.  By \eqref{ineq:pscal}, we have
\begin{equation}
 \label{ineq:u}
 2\|\uu\| \|\uy_i \wedge \uy\|
   \ge \| \langle \uu, \uy \rangle \uy_i - \langle \uu, \uy_i \rangle \uy \|
   = |\langle \uu, \uy_i \rangle| \|\uy \|
\end{equation}
for each $i\ge 0$.  Since $\|\uy_i \wedge \uy\| \sim
|\det(\uw_i)|/\|\uw_i\| \ll \|\uw_i\|^{\beta-1}$ tends to $0$ as
$i\to \infty$, we deduce from \eqref{ineq:u} that the integer
$\langle \uu, \uy_i \rangle$ must vanish for all sufficiently large
values of $i$.  This implies that $\uu=0$ because it follows from
the hypothesis $\det(\uy_0,\uy_1,\uy_2)\neq 0$ and the formula in
Proposition \ref{prop:formulesrec} (d) that any three consecutive
points of the sequence $(\uy_i)_{i\ge 0}$ are linearly independent.
Thus the coordinates of $\uy$ must be linearly independent over
$\bQ$. In particular, the first coordinate of $\uy$ is non-zero and,
dividing $\uy$ by this coordinate, we may assume that it is equal to
$1$. Then, upon denoting by $\xi$ the second coordinate of $\uy$,
the condition $\det(\uy)=0$ implies that $\uy=(1,\xi,\xi^2)$ and
thus $[\bQ(\xi):\bQ]>2$.
\end{proof}

\section{Estimates for the exponent $\omegahat_2$}
\label{section:exponents}

We first prove the following result and then deduce from it our
main theorem in \S1.

\begin{proposition}
 \label{prop:exp}
Let $(\uw_i)_{i\ge0}$ be an admissible Fibonacci sequence in $\cM$,
and let $(\uy_i)_{i\ge 0}$ be a corresponding sequence of symmetric
matrices in $\cM$.   Assume that $(\uw_i)_{i\ge0}$ is unbounded and
satisfies
\begin{equation}
 \label{ineq:norm,det}
 \|\uw_{i+1}\| \sim \|\uw_i\|^\gamma,
 \quad
 |\det(\uw_{i+1})| \sim |\det(\uw_i)|^\gamma
 \et
 \|\uw_i\|^\alpha \ll |\det(\uw_i)| \ll \|\uw_i\|^\beta
\end{equation}
for real numbers $\alpha$ and $\beta$ with $0\le \alpha\le \beta<
\gamma^{-2}$.  Assume moreover that $\tr(\uw_i)$ and $\det(\uw_i)$
are relatively prime for $i=0,1,2,3$ and that
$\det(\uy_0,\uy_1,\uy_2) \neq 0$. Then the real number $\xi$ which
comes out from the last assertion of Proposition
\ref{prop:convergence} satisfies
\begin{equation*}
 \gamma^2-\beta\gamma
 \le
 \omegahat_2(\xi)
 \le
 \gamma^2-\alpha\gamma.
\end{equation*}
\end{proposition}

\begin{proof}
Put $\uy=(1,\xi,\xi^2)$ and define the sequence $(\uz_i)_{i\ge 0}$
as in Proposition \ref{prop:formulesrec}. Since $\|\uy\|\ge 1$, the
inequality \eqref{ineq:pscal} combined with the estimates of
Proposition \ref{prop:convergence} shows that, for any point
$\uz\in\bZ^3$ and any index $i\ge 1$, we have
\begin{equation}
 \label{ineq:pscalzy}
 |\langle \uz, \uy_i \rangle|
  \le \|\uy_i\| |\langle \uz,\uy \rangle| + 2 \|\uz\| \| \uy_i\wedge \uy \|
  < c_5 \max\Big\{
       \|\uw_i\| |\langle \uz,\uy \rangle|,
       \|\uz\| \frac{|\det(\uw_i)|}{\|\uw_i\|} \Big\},
\end{equation}
with a constant $c_5>0$ which is independent of $\uz$ and $i$.
Suppose that a point $\uz\in\bZ^3$ satisfies
\begin{equation}
 \label{cond:z}
 0 < \|\uz\| \le Z_i:=c_6\|\uw_i\|
 \et
 |\langle \uz,\uy \rangle| \le \frac{|\det(\uw_{i+1})|}{\|\uw_{i+2}\|},
\end{equation}
where $c_6=c_5^{-1}|\det(\uy_2)|^{-1}$.  Using
\eqref{ineq:pscalzy} with $i$ replaced by $i+1$, we find
\begin{equation*}
 |\langle \uz, \uy_{i+1} \rangle|
  \ll |\det(\uw_i)|^\gamma \|\uw_i\|^{-1/\gamma}.
\end{equation*}
Since $|\det(\uw_i)|\ll \|\uw_i\|^\beta$ with $\beta<\gamma^{-2}$,
this gives $|\langle \uz, \uy_{i+1} \rangle|<1$ provided that $i$
is sufficiently large.  Then, the integer $\langle \uz, \uy_{i+1}
\rangle$ must be zero and, by Proposition \ref{prop:formulesrec}
(e), we deduce that $\uz=a\uz_i+b\uz_{i+1}$ for some $a,b\in\bQ$
where $b$ is given by
\begin{equation*}
 \uz_i\wedge \uz
  = b \uz_i\wedge \uz_{i+1}
  =(-1)^i b \det(\uy_0,\uy_1,\uy_2)\det(\uw_2)^{-1}\uy_{i+1}.
\end{equation*}
Since $\det(\uw_2)\uz_i\wedge \uz\in\bZ^3$ and since, by Corollary
\ref{corollary:arith} (d), the content of $\uy_{i+1}$ divides
$\det(\uy_2)/\det(\uw_2)$, this implies that $b
\det(\uy_0,\uy_1,\uy_2) \det(\uy_2)/\det(\uw_2)$ is an integer.
So, if $b$ is non-zero, it satisfies the lower bound
\begin{equation*}
 |b| \ge |\det(\uw_2)/(\det(\uy_0,\uy_1,\uy_2)\det(\uy_2))|.
\end{equation*}
We note that $\langle \uz_i, \uy_i\rangle =0$ and, by Proposition
\ref{prop:formulesrec} (d), that
\begin{equation*}
  \langle \uz_{i+1}, \uy_i \rangle
  = \frac{\det(\uy_i,\uy_{i+1},\uy_{i+2})}{\det(\uw_{i+1})}
  = (-1)^i \frac{\det(\uy_0,\uy_1,\uy_2)}{\det(\uw_2)} \det(\uw_i).
\end{equation*}
Therefore, if $b\neq 0$, the point $\uz=a\uz_i+b\uz_{i+1}$ satisfies
\begin{equation*}
 |\langle \uz, \uy_i \rangle|
  = |b| |\langle \uz_{i+1}, \uy_i \rangle|
  \ge |\det(\uy_2)|^{-1} |\det(\uw_i)|
  = c_5c_6|\det(\uw_i)|.
\end{equation*}
However, \eqref{ineq:pscalzy} and \eqref{cond:z} give
\begin{equation*}
 |\langle \uz, \uy_i \rangle|
  < c_5 \max\Big\{ \frac{|\det(\uw_{i+1})|\|\uw_i\|}{\|\uw_{i+2}\|},
                   c_6 |\det(\uw_i)| \Big\}
  = c_5c_6 |\det(\uw_i)|
\end{equation*}
if $i$ is sufficiently large, because the ratio $|\det(\uw_{i+1})|
\|\uw_i\| / \|\uw_{i+2}\| \ll \|\uw_i\|^{\beta\gamma-\gamma}$ tends
to $0$ as $i\to\infty$.   Comparison with the previous inequality
then forces $b=0$, and so we get $\uz=a\uz_i$ with $a\neq0$.  Since
$\det(\uw_2)\uz_i$ is, by Corollary \ref{corollary:arith} (e), an
integer point whose content divides
$\det(\uy_2)\det(\uy_0,\uy_1,\uy_2)$, we deduce that
$a\det(\uy_2)\det(\uy_0,\uy_1,\uy_2)/\det(\uw_2)$ is a non-zero
integer and therefore, using the second part of \eqref{ineq:conv2}
in Proposition \ref{prop:convergence}, we find that
\begin{equation*}
 |\langle \uz, \uy \rangle|
  = |a| |\langle \uz_i, \uy \rangle|
  \ge \frac{|\det(\uw_2)|}{|\det(\uy_2)\det(\uy_0,\uy_1,\uy_2)|}
      |\langle \uz_i, \uy \rangle|
  \gg \frac{|\det(\uw_{i+1})|}{\|\uw_{i+2}\|}.
\end{equation*}
Since this holds for any point $\uz$ satisfying \eqref{cond:z}
with $i$ sufficiently large, we deduce that, for any index $i\ge
0$ and any point $\uz\in\bZ^3$ with $0<\|\uz\|\le Z_i$, we have
\begin{equation*}
 |\langle \uz, \uy \rangle|
  \gg \frac{|\det(\uw_{i+1})|}{\|\uw_{i+2}\|}
  \gg \|\uw_i\|^{\gamma\alpha-\gamma^2}
  \gg Z_i^{\gamma\alpha-\gamma^2}.
\end{equation*}
This shows that $\omegahat_2(\xi)\le \gamma^2-\gamma\alpha$.

Finally, for any real number $Z\ge \|\uz_0\|$, there exists an
index $i\ge 0$ such that $\|\uz_i\| \le Z < \|\uz_{i+1}\|$ and,
for such choice of $i$, we find by Proposition
\ref{prop:convergence} that
\begin{equation*}
 |\langle \uz_i, \uy \rangle|
  \ll \frac{|\det(\uw_{i+1})|}{\|\uw_{i+2}\|}
  \ll \|\uw_i\|^{\beta\gamma-\gamma^2}
  \sim \|\uz_{i+1}\|^{\beta\gamma-\gamma^2}
  \ll Z^{\beta\gamma-\gamma^2},
\end{equation*}
showing that $\omegahat_2(\xi)\ge \gamma^2-\gamma\beta$.
\end{proof}

Let us say that a real number $\xi$ is of ``Fibonacci type'' if
there exist an unbounded Fibonacci sequence $(\uw_i)_{i\ge 0}$ in
$\cM$ and a real number $\theta$ with $\theta>1/\gamma$ such that
$\|(\xi,-1)\uw_i\|\le \|\uw_i\|^{-\theta}$ for each sufficiently
large index $i$. There are countably many such numbers, and any real
number $\xi$ obtained from Proposition \ref{prop:convergence} with
$\beta<\gamma^{-2}$ is of this type.  The following corollary shows
that the exponents $\omegahat_2(\xi)$ attached to transcendental
numbers of Fibonacci type are dense in the interval $[2,\gamma^2]$.
By Jarn\'{\i}k's formula \eqref{formula:Jarnik}, this implies our
main theorem in \S1.

\begin{corollary}
Let $t$ and $\epsilon$ be real numbers with $0<t<\gamma^{-2}$ and
$\epsilon>0$.  Then, there exist a transcendental real number
$\xi$ and an unbounded Fibonacci sequence $(\uw_i)_{i\ge 0}$ in
$\cM$ which satisfy
\begin{itemize}
 \item[(a)] $\|(\xi,-1)\uw_i\|\le \|\uw_i\|^{-1+t}$ for each
  sufficiently large $i$,
 \item[(b)] $\gamma^2-t\gamma \le \omegahat_2(\xi) \le
              \gamma^2-(t-\epsilon)\gamma$.
\end{itemize}
\end{corollary}

\begin{proof}
Since $t<1$, there exist integers $k$ and $\ell$ with $0<\ell<k$ and
$t-\epsilon \le \ell/(k+2)\le \ell/k < t$. For such a choice of $k$
and $\ell$, consider the Fibonacci sequence $(\uw_i)_{i\ge 0}$ of
Example \ref{example1} with parameters $a=2^\ell$, $b=2^{k-\ell}-1$
and $c=2^{k-\ell}$. According to Example \ref{example2}, $\uw_i$ has
relatively prime trace and determinant for each $i\ge 0$ and the
corresponding sequence of symmetric matrices $(\uy_i)_{i\ge 0}$
satisfies $\det(\uy_0,\uy_1,\uy_2)=2^{4\ell}\neq 0$.  Moreover,
Example \ref{example3} shows that $(\uw_i)_{i\ge 0}$ is unbounded
and satisfies the estimates \eqref{ineq:norm,det} of Proposition
\ref{prop:exp} with $\alpha=\ell/(k+2)$ and $\beta=\ell/k$ (note
that the example provides a slightly larger value for $\alpha$). So,
Proposition \ref{prop:exp} applies and shows that the corresponding
real number $\xi$ constructed by Proposition \ref{prop:convergence}
satisfies the above condition (b). In particular, $\xi$ is
transcendental since $\omegahat_2(\xi)>2$. Moreover, since
$\|(\xi,-1)\uw_i\| \sim \|(\xi,-1)\uy_i\| \sim \|\uy_i\wedge\uy\|$,
the first estimate in \eqref{ineq:conv2} leads to (a).
\end{proof}


\end{document}